\documentclass[11pt]{article}

\usepackage{amscd,amsmath, amssymb, fancyhdr, epsfig}

\numberwithin{equation}{section}

\newcommand{\version}{version 2.0,\ \   14.02.2008}

\def\eqref#1{(\ref{#1})}
\newcommand{\goth}{\mathfrak}

\newcommand{\arrow}{{\:\longrightarrow\:}}

\newcommand{\C}{{\Bbb C}}
\newcommand{\R}{{\Bbb R}}

\newcommand{\6}{\partial}
\def\1{\sqrt{-1}\:}
\newcommand{\restrict}[1]{{\left|_{{\phantom{|}\!\!}_{#1}}\right.}}
\newcommand{\cntrct}                
{\hspace{2pt}\raisebox{1pt}{\text{$\lrcorner$}}\hspace{2pt}}

\renewcommand{\c}[1]{{\cal #1}}
\newcommand{\calo}{{\cal O}}
\newcommand{\cac}{{\cal C}}

\renewcommand{\tilde}{\widetilde}
\renewcommand{\bar}{\overline}
\renewcommand{\phi}{\varphi}
\renewcommand{\epsilon}{\varepsilon}
\renewcommand{\geq}{\geqslant}
\renewcommand{\leq}{\leqslant}
\renewcommand{\max}{{\rm max}}

\newcommand{\Vol}{\operatorname{Vol}}

\newcommand{\Sing}{\operatorname{Sing}}

\newcommand{\dist}{\operatorname{dist}}
\newcommand{\diag}{\operatorname{\sf diag}}

\newcounter{Mycounter}[section]
\newcounter{lemma}[section]
\renewcommand{\thelemma}{{Lemma \thesection.\arabic{lemma}}}
\newcommand{\lemma}{%
    \setcounter{lemma}{\value{Mycounter}}
    \refstepcounter{lemma}
    \stepcounter{Mycounter}
    {\noindent \bf \thelemma:\ }}

\newcounter{claim}[section]
\renewcommand{\theclaim}{{Claim \thesection.\arabic{claim}}}
\newcommand{\claim}{%
    \setcounter{claim}{\value{Mycounter}}
    \refstepcounter{claim}
    \stepcounter{Mycounter}
    {\noindent \bf \theclaim:\ }}

\newcounter{sublemma}[section]
\renewcommand{\thesublemma}{{Sublemma \thesection.\arabic{sublemma}}}
\newcommand{\sublemma}{%
    \setcounter{sublemma}{\value{Mycounter}}
    \refstepcounter{sublemma}
    \stepcounter{Mycounter}
    {\noindent \bf \thesublemma:\ }}

\newcounter{corollary}[section]
\renewcommand{\thecorollary}{{Corollary \thesection.\arabic{corollary}}}
\newcommand{\corollary}{%
    \setcounter{corollary}{\value{Mycounter}}
    \refstepcounter{corollary}
    \stepcounter{Mycounter}
    {\noindent \bf \thecorollary:\ }}

\newcounter{theorem}[section]
\renewcommand{\thetheorem}{{Theorem \thesection.\arabic{theorem}}}
\newcommand{\theorem}{%
    \setcounter{theorem}{\value{Mycounter}}
    \refstepcounter{theorem}
    \stepcounter{Mycounter}
    {\noindent \bf \thetheorem:\ }}

\newcounter{conjecture}[section]
\newcounter{proposition}[section]
\renewcommand{\theproposition}
      {{Proposition \thesection.\arabic{proposition}}}
\newcommand{\proposition}{%
    \setcounter{proposition}{\value{Mycounter}}
    \refstepcounter{proposition}
    \stepcounter{Mycounter}
    {\noindent \bf \theproposition:\ }}

\newcounter{definition}[section]
\renewcommand{\thedefinition}
      {{Definition~\thesection.\arabic{definition}}}
\newcommand{\definition}{%
    \setcounter{definition}{\value{Mycounter}}
    \refstepcounter{definition}
    \stepcounter{Mycounter}
    {\noindent \bf \thedefinition:\ }}

\newcounter{example}[section]
\newcounter{remark}[section]
\renewcommand{\theremark}{{Remark \thesection.\arabic{remark}}}
\newcommand{\remark}{%
    \setcounter{remark}{\value{Mycounter}}
    \refstepcounter{remark}
    \stepcounter{Mycounter}
    {\noindent \bf \theremark:\ }}

\newcounter{problem}[section]
\newcounter{question}[section]
\makeatletter

\@addtoreset{equation}{section} \@addtoreset{footnote}{section}
\makeatother

\def\blacksquare{\hbox{\vrule width 5pt height 5pt depth 0pt}}
\def\endproof{\blacksquare}

\begin{document}
\begin{center}
{\LARGE\bf
Plurisubharmonic functions in calibrated\\[1mm] geometry and $q$-convexity\\[4mm]
}

 Misha
Verbitsky\footnote{Misha Verbitsky is 
supported by CRDF grant RM1-2354-MO02.}

\end{center}

{\small \hspace{0.15\linewidth}
\begin{minipage}[t]{0.7\linewidth}
{\bf Abstract} \\ 
Let $(M, \omega)$ be a K\"ahler manifold. 
An integrable function $\phi$ on $M$ is called
$\omega^q$-plurisubharmonic if the current
$dd^c\phi\wedge \omega^{q-1}$ is positive.
We prove that $\phi$ is $\omega^q$-plurisubharmonic
if and only if $\phi$ is subharmonic on all
$q$-dimensional complex subvarieties. We prove
that a $\omega^q$-plurisubharmonic function
is $q$-convex, and admits a local approximation
by smooth, $\omega^q$-plurisubharmonic functions.
For any closed subvariety $Z\subset M$, 
$\dim_\C Z\leq q-1$, there exists a strictly
$\omega^q$-plurisubharmonic function in a neighbourhood
of $Z$ (this result is known for $q$-convex
functions). This theorem is used to give a new
proof of Sibony's lemma on integrability of
positive closed $(p,p)$-forms which are integrable
outside of a complex subvariety of codimension  
$\geq p+1$.
\end{minipage}
}

\tableofcontents


\section{Introduction}


\subsection{Positive forms and $\omega^q$-plurisubharmonic functions}

The notion of $\omega^q$-plurisubharmonic function has
many facets, and these functions can be defined in many
different ways. The main reason this notion is considered
comes from the theory of plurisubharmonic functions
on calibrated manifolds (see Subsection
\ref{_psh_calibra_Subsection_}). However,
$\omega^q$-plurisubharmonic functions are
very useful even outside of the theory of 
calibrations.

Let $(M, \omega)$ be a K\"ahler manifold,
and $\phi:\; M \arrow \R$ a smooth function.
Consider the form $dd^c\phi \in \Lambda^{1,1}(M)$,
\[
dd^c\phi = 2 \1 \6\bar\6 \phi.
\]
Using the standard linear-algebraic argument,
in a contractible open set $U$ we might find an orthonormal frame
$ \xi_1, ... \xi_n \in \Lambda^{1,0}(U)$, such that
\begin{equation}\label{_dd^c_locally_ortho_Equation_}
dd^c\phi = -\1\sum_i \alpha_i \xi_i\wedge\bar\xi_i.
\end{equation}
Here, $\alpha_i$ are real functions,
and the set $\{\alpha_i\}$ is independent from the
choice of a frame.
A function $\phi$ is called {\bf $\omega^q$-plurisubharmonic}
if the sum of any $q$ eigenvalues is positive:
\begin{equation}\label{_sum_smallest_Equation_}
\sum_{k=1}^q \alpha_{i_k} \geq 0,
\end{equation}
for any $k$-tuple $i_1 < i_2 < ... < i_q$. 
This condition implies that at most $q-1$ eigenvalues
of $dd^c\phi$ can be negative. Such functions are known
as {\bf $q$-convex}. Indeed,
$\omega^q$-plurisubharmonicity
is in many aspects similar to $q$-convexity.

For $q=1$, being $\omega^q$-plurisubharmonic
is equivalent to being plurisubharmonic.

Recall that a real $(p,p)$-form $\eta$
on a complex manifold is called {\bf weakly positive}
if for any complex subspace $V\subset T_c M$, 
$\dim_\C V=p$, the restriction $\rho\restrict V$
is a non-negative volume form. Equivalently,
this means that 
\[ 
  (\1)^p\rho(x_1, \bar x_1, x_2, \bar x_2, ... x_p, \bar
  x_p)\geq 0,
\]
for any vectors $x_1, ... x_p\in T_x^{1,0}M$.
A form is called {\bf strongly positive} if it can 
be expressed as a sum
\[
\eta = (-\1)^p\sum_{i_1, ... i_p} 
\alpha_{i_1, ... i_p} \xi_{i_1} \wedge \bar\xi_{i_1}\wedge ... 
\wedge \xi_{i_p} \wedge \bar\xi_{i_p}, \ \  
\]
running over some set of $p$-tuples 
$\xi_{i_1}, \xi_{i_2}, ..., \xi_{i_p}\in \Lambda^{1,0}(M)$,
with $\alpha_{i_1, ... i_p}$ real and non-negative functions on $M$.

The strongly positive and the weakly positive forms
form closed, convex cones in the space 
$\Lambda^{p,p}(M,\R)$ of real $(p,p)$-forms.
These two cones are dual with respect to the Poincare pairing
\[
\Lambda^{p,p}(M,\R) \times \Lambda^{n-p,n-p}(M,\R)\arrow \Lambda^{n,n}(M,\R)
\]
where $n=\dim_\C M$.
For (1,1)-forms and $(n-1,n-1)$-forms,
the strong positivity is also equivalent
to weak positivity.

Throughout this paper, we are mostly interested in 
differential forms of type $\eta = \nu \wedge \omega^k$,
where $\nu\in \Lambda^{1,1}(M)$ is a real (1,1)-form.
For such forms, strong positivity is equivalent
to weak positivity.

\hfill

\claim \label{_strong_posi_equi_weak_Claim_}
Let $(M, \omega)$ be a K\"ahler manifold, and 
$\nu\in \Lambda^{1,1}(M)$ a real (1,1)-form.
Then, for any $k\geq 0$, $\nu\wedge \omega^k$
is weakly positive if and only if it is strongly positive.

\hfill

{\bf Proof:} Consider a decomposition similar to 
\eqref{_dd^c_locally_ortho_Equation_},
\[\nu = -\1\sum_i \alpha_i \xi_i\wedge\bar\xi_i.
\]
Then $\nu\wedge \omega^k$ can be written as
\begin{equation}\label{_nu_wedge_omega^k_Equation_}
\nu\wedge \omega^k = (\1)^k \sum_{i_1< i_2 < ... <
  i_{k+1} } 
\left(\sum_{j=1}^{k+1} \alpha_{i_j}\right) 
\xi_{i_1}\wedge\bar\xi_{i_1}\wedge\xi_{i_2}\wedge\bar\xi_{i_2}
\wedge ...\wedge \xi_{i_{k+1}}\wedge\bar\xi_{i_{k+1}},
\end{equation}
where the first sum is taken over all $(k+1)$-tuples.

For such a form, weak positivity is equivalent to strong
positivity and equivalent to positivity of all
coefficients $\sum_{j=1}^{k+1} \alpha_{i_j}$. \endproof

\hfill

From \eqref{_nu_wedge_omega^k_Equation_}
it is clear that $\phi$ is $\omega^q$-plurisubharmonic
if and only if $dd^c\phi\wedge \omega^{q-1}$
is positive (see \ref{_omega^q_psh_via_eigenva_Theorem_}
for a detailed argument). 

In \ref{_subharmo_and_omega^q_psh_Theorem_},
we prove that a function $\phi$ is  $\omega^q$-plurisubharmonic
if and only if $\phi$ is subharmonic on all 
germs of $q$-dimensional complex subvarieties
of $M$.

For a domain in $\C^n$, the notion of 
$\omega^q$-plurisubharmonicity was introduced and studied
by Z. Khusanov, under the name ``$q$-subharmonic'' 
(\cite{_Khusanov_1_}, \cite{_Khusanov_2_},
\cite{_Khusanov_3_}). These functions were considered
by Z. B\l ocki in \cite{_Blocki_}, also for domains in
$\C^n$, and used in connection with a generalized
version of Monge-Ampere equation, called ``complex 
Hessian equation''.

\subsection{Strictly positive forms}
\label{_strict_posi_Subsection_}

Let $(M, \omega)$ be a K\"ahler manifold.
The sets of strongly and weakly positive $(p,p)$-forms
form cones $\cac_s$ and $\cac_w$ in $\Lambda^{p,p}(M)$.
We consider $\cac_s$, $\cac_w$ as topological spaces
with open-compact topology. A form $\eta$
lies in the interior part of $\cac_s$, $\cac_w$,
if for any compact $K\subset M$ there exists
$\epsilon >0$ such that $(\eta- \epsilon \omega^p)\restrict K$
lies in $\cac_s$, $\cac_w$. Such a form is called 
{\bf strictly positive} (strictly strongly
positive and strictly weakly positive).

Exhausting $M$ with compact sets,
we obtain that $\eta$ is strictly positive if
and only if there exists a continuous positive
function $h:\; M \arrow ]0, \infty[$ such that
$\eta-  h\omega^p$ is positive. 

A function $\phi$ is called {\bf strictly 
$\omega^q$-plurisubharmonic} if $dd^c\phi\wedge \omega^{q-1}$
is strictly positive (strong positivity and weak
positivity are equivalent in this situation: see 
\ref{_strong_posi_equi_weak_Claim_}).
This is equivalent to the strictness of
the inequality \eqref{_sum_smallest_Equation_}.
For $q=1$, this gives a definition of strictly
plurisubharmonic functions.

\subsection{Continuous $\omega^q$-plurisubharmonic
functions}

For any locally integrable form $\eta$, 
$dd^c\eta$ is a well defined current. An 
integrable function $f$ is called plurisubharmonic
if $dd^cf$ is a positive current. Continuous
strictly plurisubharmonic functions can be approximated
by smooth strictly plurisubharmonic functions. This
approximation is performed by taking coordinates
and smoothing a function by a convolution with 
a smoothing kernel. 

For $\omega^q$-plurisubharmonic
functions, $q>0$, this argument fails. The smoothing
operation described above is a convolution with a 
smooth kernel on $G\times M$, where $G$ is a group 
of local automorphisms of $M$. For such a
convolution to be smooth, this group should
act transitively in a neghbourhood of a point.
However, for $q>0$ a notion of $\omega^q$-plurisubharmonicity
depends on the choice of a K\"ahler form,
and for $\omega$ sufficiently general, 
$(M, \omega)$ has no local isometries.

A similar problem arises when one tries to approximate
continuous subharmonic functions on arbitrary
Riemannian manifold by smooth subharmonic functions.
Recall that a $L^1$-integrable function
$f$ is {\bf subharmonic} if the current
$\Delta f$ is positive. On a K\"ahler manifold $(M, \omega)$,
$\dim_\C M=n$, a function is $\omega^n$-plurisubharmonic
if and only if it is subharmonic; this is clear,
because $\Delta f \wedge \omega^n= dd^cf \wedge \omega^{n-1}$.

Green and Wu (\cite{_Green_Wu:approx_}) used the heat equation
to prodice smooth approximations for all continuous
subharmonic functions.  We modify their argument, obtaining
a smooth approximation for all strictly 
$\omega^q$-plu\-ri\-sub\-har\-monic functions 
(\ref{_local_appro_Theorem_},
\ref{_Global_appro_omega_psh_Theorem_}).

\subsection{Strictly $\omega^q$-plurisubharmonic functions
in a neighbourhood of a subvariety}

One of utilities of $q$-convex functions comes from the following
theorem.

\hfill

\theorem
Let $M$ be a complex manifold, and $Z\subset M$ 
a closed complex subvariety, $\dim_\C Z\leq q$.
Then there exists an open neighbourhood 
$V\supset Z$ and a smooth strictly 
$q$-convex function $\phi:\; V \arrow R$.

\hfill

{\bf Proof:} \cite{_Barlet_}, \cite{_Demailly:q-convex_}. \endproof

\hfill

However, $q$-convex functions are not very convenient to use.
A sum, maximum, or regularized maximum of two $q$-convex 
functions  is no longer $q$-convex. Also, continuous
$q$-convex function do not admit smooth approximations.

On contrast, a sum, maximum, and regularized maximum
of $\omega^q$-plu\-ri\-sub\-harmonic functions is again 
$\omega^q$-plu\-ri\-sub\-har\-monic; smooth approximation is
 also possible. 

Existence of strictly $q$-convex exhaustion functions leads
to well-known restrictions on topology and geometry of a manifold $M$.
Using the Morse theory, such a function gives a cell decomposition
of this manifold. This can be used to show that $M$ is homotopy equivalent
to a CW-complex of dimension at most $\dim_\C M +q -1$.
Even more information can be obtained using 
vanishing theorems obtained from $L^2$-estimates
(see e.g. \cite{_Demailly:q-convex_}, \cite{_Ohsawa_}).

For $q>1$, $\omega^q$-plurisubharmonicity is obviuosly 
stronger than $q$-convexity. However, it seems that there
are no additional geometric restrictions on existence
of strictly $\omega^q$-plurisubharmonic functions. At least,
a version of \ref{_omegaq-psh-in-neigh_Theorem_} is 
valid in this case.

\hfill

\theorem
Let $Z\subset M$ be a closed complex subvariety of a 
K\"ahler manifold. Then there exists a strictly $\omega^q$-plurisubharmonic 
function $\phi$ on an open neighbourhood of $Z$ in $M$.

\hfill

{\bf Proof:} This is \ref{_omega^q_psh_in_neigh_Proposition_}.
\endproof 

\hfill

When $Z$ is compact, $\phi$ can be chosen exhausting,
in appropriate neighbourhood of $Z$ (\ref{_omegaq-psh-in-neigh_Theorem_}).

\subsection{Similar results in literature}

After this paper was finished, Mihnea Col\c toiu
sent me the reference to a number of papers
were similar results were discussed.

The strictly $\omega^q$-plurisubharmonic functions
were studied by H. Wu, in \cite{_Wu:q_complete_},
under the name of {\em class $\Psi(q)$-functions.}
Wu has shown that strictly $\omega^q$-plurisubharmonic functions
admit smooth approximations (\ref{_Global_appro_omega_psh_Theorem_}).

In \cite{_Napier_Ramachandran_}, T. Napier and M. Ramachandran
applied results of Wu to prove Bochner-Hartogs dichotomy
for certain classes of K\"ahler manifolds, showing that a Kaehler
manifold which admits a continuous exhausting plurisubharmonic 
function and has exactly one end either admits a holomorphic
map to a Riemannian surface, or satisfies $H^1_c(M, \calo)=0$,
where $H^1_c(M, \calo)$ is a group of holomorphic cohomology
with compact support. Napier and Ramachandran consider
$\omega^q$-plurisubharmonic functions, under the name of
{\em $q$-plurisubharmonic functions}. They state
a version of \ref{_omegaq-psh-in-neigh_Theorem_}
and sketch its proof, referring to Demailly and Col\c toiu, who proved a
similar result for $q$-convex functions in 
\cite{_Demailly:q-convex_} and \cite{_Coltoiu_}.


\section{$\omega^q$-plurisubharmonic functions}


\subsection{Plurisubharmonic functions in calibrated
  geometry}
\label{_psh_calibra_Subsection_}

Since their introduction by 
Harvey and Lawson in \cite{_Harvey_Lawson:Calibrated_},
calibrated geometries became a mainstay of 
much of modern differential geometry and
some of string physics. Through calibrations,
the arguments of K\"ahler geometry can be applied
to manifolds with special holonomy, and
 other interesting differential-geometric structures.

\hfill

\definition
Let $M$ be a Riemannian manifold. A {\bf calibration} on
$M$ is a closed real $k$-form $\rho\in \Lambda^k(M)$
such that for any $k$-dimensional subspace $V\subset T_xM$, 
we have 
\begin{equation}\label{_calibra_subspa_Equation_}
\bigg|\rho|_ V\bigg| \leq |\Vol_V|,
\end{equation}
where $\Vol_V$ is a Riemannian volume form on $V$.

A space $V\subset T_xM$ is called {\bf calibrated by
$\rho$} is $\bigg|\rho|_V\bigg| =|\Vol_V|$.
A $k$-dimensional submanifold $Z\subset M$ is 
{\bf calibrated by $\rho$} if $T_zZ\subset T_zM$ 
is calibrated, for all $z\in Z$.

Examples of well-known calibrations include 
the K\"ahler form and its powers (normalized as
$\frac {1}{q!}\omega^q$); the fundamental 3-form and
4-form on a $G_2$-manifold; real part of the holomorphic
volume form on a Calabi-Yau manifold; and so on. 

For a detailed study of calibrations on manifolds,
the reader can consult the books \cite{_Harvey:Spinors_}
and \cite{_Joyce:Calibrated_}.

The main utility of calibrations is due to the
following important theorem.

\hfill

\theorem 
(\cite{_Harvey_Lawson:Calibrated_})
Let $(M, g, \rho)$ be a Riemannian manifold equipped with
a calibration, and $Z\subset M$ an oriented submanifold. Then
\[
\int_Z \Vol_Z \geq \bigg| \int_Z\rho\bigg|.
\]
Moreover, the equality happens only when $Z$ is 
calibrated by $\rho$.

\hfill

{\bf Proof:} Follows immediately from
\eqref{_calibra_subspa_Equation_}.
\endproof

\hfill

This result is remarkable, because for $Z$ compact,
the integral $\int_Z\rho$ is a cohomological invariant. Therefore,
compact calibrated manifolds are minimal, and, moreover,
minimize the Riemannian volume in their cohomology class.

When $\rho = \frac 1 {q!} \omega^q$, where 
$\omega$ is a K\"ahler form, a manifold $Z\subset M$
is calibrated with respect to $\rho$ if and only if $Z$ is
complex analytic. 

\hfill

In \cite{_Harvey_Lawson:Psh_}, Harvey and Lawson
defined the notion of a plurisubharmonic function
associated with a calibration. In the present paper, we are
interested in parallel calibrations, that is, calibrations
satisfying $\nabla \rho=0$, where $\nabla$ is the
Levi-Civita connection. In \cite{_Harvey_Lawson:Psh_}, 
plurisubharmonic functions are defined with respect to
an arbitrary calibration, but for parallel ones the
definition is simpler. We restrict ourselves to
parallel calibrations.

Consider a parallel calibration $\rho \in \Lambda^k(M)$,
and let $\Lambda^1(M) \stackrel C \arrow \Lambda^{k-1}(M)$
be the map $\theta \arrow \rho \cntrct \theta^\sharp$,
where $\cntrct$ is a contraction, and $\theta^\sharp$
the dual vector field. When $\rho$ is the K\"ahler form,
this is a complex structure operator; when $\rho$ is
the fundamental 3-form on a $G_2$-manifold, $C$ is
dual to the vector product (\cite{_Verbitsky:G2_}).  

Consider an operator $d_c:= \{ d, C\}$, where
$\{\ \}$ denotes the supercommutator, 
\[
\{ d, C\} = dC - (-1)^{\deg \rho} Cd.
\]
The operator $dd_c:\; C^\infty M \arrow \Lambda^k M$ 
(introduced independently in \cite{_Verbitsky:G2_})
is a natural analogue of the operator $dd^c$ in complex
geometry. When $\rho$ is a K\"ahler form,
$dd_c$ is equal to $dd^c$.

Given a calibrated subspace $V\subset T_xM$,
$\rho\restrict V$ is non-degenerate, and gives an
orientation. A smooth function $\phi:\; M \arrow \R$
is called {\bf plurisubharmonic with respect to $\rho$},
if $dd_c \phi$ is positive on any calibrated subspace
$V \subset T_x M$. We extend this definition to
integrable functions in 
Subsection \ref{_conti_omega^q_psh_Subsection_}. When
the restriction $dd_c\phi\restrict V$ is non-degenerate on all
$\rho$-calibrated subspaces, $\phi$ is called
{\bf strictly plurisubharmonic with respect to $\rho$}.

We are going to explore this notion in the case
when $\rho = \frac 1 {q!} \omega^q$, where $\omega$ is a 
K\"ahler form. 

\subsection{$\omega^q$-plurisubharmonic functions and
  $q$-convexity}

From now on, $(M, \omega)$ is a K\"ahler manifold,
considered as a calibrated manifold with the calibration 
$\rho = \frac 1 {q!} \omega^q$.
A function is called {\bf $\omega^q$-\-plu\-ri\-sub\-har\-mo\-nic}
if it is plurisubharmonic with respect to the calibration
$\frac 1 {q!} \omega^q$. 

\hfill

\claim\label{_d_c_via_d^c_Claim_}
In these assumptions, 
\[ 
d_c(\phi) = \frac1 {(q-1)!} d^c\phi \wedge \omega^{q-1},
\]
where $d^c = I d I^{-1} $ is the usual twisted de Rham
differential.

\hfill

{\bf Proof:} Clearly, 
\[ d^c\phi= \omega \cntrct (d\phi)^\sharp, \]
and 
\[  d_c\phi = q \frac 1 {q!} \omega^{q-1} \wedge (\omega \cntrct
(d\phi)^\sharp).
\]
From these two equations, \ref{_d_c_via_d^c_Claim_} 
clearly follows. \endproof

\hfill

From \ref{_d_c_via_d^c_Claim_}, we obtain that the
$\omega^q$-plurisubharmonicity of a function $\phi$
depends on the pseudo-Hermitian $(1,1)$-form $dd^c\phi$.
In some orthonormal basis $\xi_1, ... \xi_n$ 
in $\Lambda^{1,0}_z(M)$, this form can be written 
as $dd^c \phi = \1\sum_i \alpha_i \xi_i\wedge \bar\xi_i$,
with $\alpha_i$ real numbers. We call these numbers
{\bf the eigenvalues of $dd^c \phi$} at $z\in M$. Clearly, the set
$\{\alpha_i\}$ is independent from the choice of the
orthonormal basis $\xi_i$.

\hfill

\theorem\label{_omega^q_psh_via_eigenva_Theorem_}
Let $(M, \omega)$ be a K\"ahler manifold, and 
$\phi:\; M \arrow \R$ a smooth function. 
Then the following conditions are equivalent.
\begin{description}
\item[(i)] $\phi$ is $\omega^q$-plurisubharmonic
\item[(ii)] $dd^c\phi\wedge \omega^{q-1}$ is weakly
positive
\item[(ii)] $dd^c\phi\wedge \omega^{q-1}$ is strongly
positive
\item[(iv)] Suppose we order the eigenvalues of
  $dd^c\phi$ as follows: 
$\alpha_1 \leq \alpha_2 \leq  ... \leq \alpha_n$.
Then
\begin{equation}\label{_omega^q-psh_eigen_Equation_}
\sum_{i=1}^{q}\alpha_i>0.
\end{equation}
\end{description}

\remark
The inequality \eqref{_omega^q-psh_eigen_Equation_}
implies that $\alpha_q, \alpha_{q+1}, ... \alpha_n$
are non-negative. Indeed, if $\alpha_q$ is negative, then
$\alpha_1, ..., \alpha_{q-1}$ are also negative, but
then $\sum_{i=1}^{q}\alpha_i<0$.

\hfill

{\bf Proof of \ref{_omega^q_psh_via_eigenva_Theorem_}:}
The equivalence (i) $\Leftrightarrow$ (ii) is clear from
\ref{_d_c_via_d^c_Claim_}. Indeed, a plane $V\subset T_x
M$ is $\omega^q$-calibrated if and only if it is complex;
therefore, $\phi$ is $\omega^q$-plurisubharmonic if
and only if $dd_c \phi= \frac 1 {(q-1)!} dd^c\phi \wedge \omega^{q-1}$
is weakly positive. The equivalence (ii) $\Leftrightarrow$ (iii)
follows from \ref{_strong_posi_equi_weak_Claim_}. Now, let $\xi_1, ... \xi_n$ 
be an orthonormal basis in $\Lambda^{1,0}(M)$, such that
$dd^c \phi = \1\sum_i \alpha_i \xi_i\wedge \bar\xi_i$.
Then
\begin{equation}\label{_dd_c_phi_coord_Equation_}
 \frac 1 {(q-1)!} dd^c\phi \wedge \omega^{q-1} = 
 \sum_{i_1< i_2 < ... < i_q} \bigg( \sum_{k=1}^q \alpha_{i_k}\bigg)
 (\1)^q
\xi_{i_1}\wedge \bar\xi_{i_1} \wedge ... \wedge \xi_{i_q}\wedge \bar\xi_{i_q}.
\end{equation}
Clearly, when all the coefficients $\sum_{k=1}^q \alpha_{i_k}$
are positive, the form \eqref{_dd_c_phi_coord_Equation_}
is also positive, being a sum of positive monomials. 
The converse is also true.
Indeed, 
\[ 
 (\1)^q \frac 1 {(q-1)!} dd^c\phi \wedge \omega^{q-1}
 (\xi_{i_1}^\sharp, \bar\xi_{i_1}^\sharp,
 \xi_{i_2}^\sharp, \bar\xi_{i_2}^\sharp, ..., 
\xi_{i_q}^\sharp, \bar\xi_{i_q}^\sharp) = \sum_{k=1}^q \alpha_{i_k},
\]
where $\xi_{i_2}^\sharp\in TM$ is the basis dual to $\xi_i$.
Therefore, if one of the coefficients
$\sum_{k=1}^q \alpha_{i_k}$ is negative, then
$\frac 1 {(q-1)!} dd^c\phi \wedge \omega^{q-1}$ is
negative on the space 
$V=\langle \xi_{i_1}^\sharp, \xi_{i_2}^\sharp, ..., \xi_{i_q}^\sharp\rangle$
generated by the corresponding $q$-tuple.

We obtained that the positivity of the form 
$\frac 1 {(q-1)!} dd^c\phi \wedge \omega^{q-1}$
is equivalent to non-negativity of all the coefficients
$\sum_{k=1}^q \alpha_{i_k}$. The smallest of these
coefficients is $\sum_{k=1}^q \alpha_{k}$.
Therefore, positivity of 
$\frac 1 {(q-1)!} dd^c\phi \wedge \omega^{q-1}$
is equivalent to $\sum_{k=1}^q \alpha_{k}\geq 0$.
We proved \ref{_omega^q_psh_via_eigenva_Theorem_}.
\endproof

\hfill

\remark\label{_strict_omega^q_psh_eigenva_Remark_}
In these assumptions,
$\phi$ is strictly $\omega^q$-plurisubharmonic
if and only if the sum of $q$ smallest
eigenvalues $\sum_{k=1}^q \alpha_{k}$ is positive everywhere; 
this is proven by the same argument as the above theorem.

\hfill

In complex geometry, the notion of $q$-convexity is 
used quite often. A function $\phi:\; M \arrow \R$
on a complex manifold is called {\bf strongly $q$-convex}
(see e.g. \cite{_Demailly:q-convex_}) if $dd^c\phi$ has at
most $(q-1)$ negative or zero eigenvalues. 
Even if the eigenvalues of $dd^c\phi$ 
are defined in terms of a K\"ahler (or a Hermitian) form,
their sign is independent from this choice. Therefore,
the $q$-convexity of a function depends only on the
complex structure on $M$.

From \ref{_strict_omega^q_psh_eigenva_Remark_},
it is clear that strictly $\omega^q$-plurisubharmonic
functions are strongly $q$-convex. 
The converse is not true: $q$-convexity 
does not necessarily imply the inequality 
$\sum_{k=1}^q \alpha_{k}>0$ for the
sum of $q$ smallest eigenvalues of $dd^c\phi$.

However, unlike plurisubharmonic functions, 
$q$-convex functions are not very convenient to use.
A sum of two $q$-convex functions is
no longer $q$-convex. Moreover, continuous $q$-convex
functions (defined as a function which can be locally
expressed as a maximum of a funite number of smooth
$q$-convex functions)\footnote{Such functions are also
known as {\bf $q$-convex functions with corners}.}
do not necessarily have a smooth $q$-convex approximation
(\cite{_Dieder_Fornae:smoothing_}, \cite{_M_Peternell:contin_q_con_}).
All of these problems are rectified if we consider 
$\omega^q$-plurisubharmonic functions instead of
$q$-convex: sum and maximum
of $\omega^q$-plurisubharmonic is again
$\omega^q$-plurisubharmonic, and any strict
$\omega^q$-\-plu\-ri\-sub\-har\-mo\-nic function admits a smooth
approximation (see Section \ref{_omega_psh_properties_Section_}).



\section{Properties of $\omega^q$-plurisubharmonic
  functions}
\label{_omega_psh_properties_Section_}

\subsection{Continuous $\omega^q$-plurisubharmonic
  functions}
\label{_conti_omega^q_psh_Subsection_}

Just like it happens in the case of the usual plurisubharmonic
functions, the definition of $\omega^q$-plurisubharmonic
functions can be generalized, to include continuous $\omega^q$-plurisubharmonic
functions. 

The following claim is trivial.

\hfill

\claim\label{_psh_via_int_Claim_}
Let $M$ be a K\"ahler manifold.
A smooth function $\phi$ is $\omega^q$-plurisubharmonic
if and only if 
\begin{equation}\label{_dd_c_Stokes_Equation_}
\int_M \phi \omega^{q-1} \wedge dd^c \alpha \geq 0,
\end{equation}
for any strongly positive form $\alpha$ with compact
support.

\hfill

{\bf Proof:} We use
\[
\int_M  \phi \omega^{q-1} \wedge dd^c \alpha = 
\int_M  dd^c\phi \wedge \omega^{q-1} \wedge  \alpha,
\]
which follows from Stokes' theorem. \endproof

\hfill

Now, \eqref{_dd_c_Stokes_Equation_}
can be used to define the $\omega^q$-plurisubharmonicity
for continuous functions (see also \cite{_Harvey_Lawson:Dua_}).

\hfill

\definition
Let $(M, \omega)$ be a K\"ahler manifold,
and $\phi:\; M \arrow \R$ a locally integrable
function. Then $\phi$ is called {\bf
  $\omega^q$-plurisubharmonic} if 
\[
\int_M \phi \omega^{q-1} \wedge dd^c \alpha \geq 0
\]
for any positive form $\alpha$ with compact support.

From \ref{_psh_via_int_Claim_}, 
it is obvious that for smooth $\phi$ this definition
agrees with the one we gave previously. Also, the following
statement is apparent.

\hfill

\claim\label{_limit_psh_psh_Claim_}
Let $\{\phi_i\}$ be a sequence of
$\omega^q$-plurisubharmonic functions which locally
converge to $\phi$ in $L^1$-topology. Then $\phi$
is also $\omega^q$-plurisubharmonic.

\endproof

\subsection{Local approximation for
  $\omega^q$-plurisubharmonic functions}
\label{_local_appro_Subsection_}

The $\omega^q$-plurisubharmonic functions admit a smooth
approximation, both locally 
(\ref{_local_appro_Theorem_}) and 
globally (\ref{_Global_appro_omega_psh_Theorem_}).
The arguments we use are taken from \cite{_Green_Wu:approx_},
where the same result is proven for subharmonic functions
on Riemannian manifolds.

\hfill

\theorem\label{_local_appro_Theorem_}
Let $(M, \omega)$ be an open K\"ahler manifold, biholomorphic to
an open ball, $K\subset M$ a compact subset, and
$\phi:\; M \arrow \R$ a continuous
$\omega^q$-plurisubharmonic function. Then there
exists a sequence $\{\phi_i\}$ of smooth functions
on $M$, $\omega^q$-plurisubharmonic in some neighbourhood of $K$,
and uniformly converging to $\phi$ on $K$.

\hfill

{\bf Proof:} The proof of \ref{_local_appro_Theorem_}
uses the notion of a {\em heat kernel}. We recall the
necessary results, which are well known (see e.g.
\cite{_Hain:heat_kernel_}). 

\hfill

\theorem\label{_Heat_Riemannian_Theorem_}
Let $M$ be a  compact Riemannian manifold,
$\dim_\R M = n$, and $\eta \in \Lambda^k(M)$ an
$L^1$-integrable $k$-form. Then the equation
\[
\frac{d\eta_t}{dt}= -\Delta \eta_t, \ \  \eta_0 = \eta,
\]
has a unique solution $\eta_t \in \Lambda^k(M)$,
smoothly depending on $t\in \R^{\geq 0}$. Moreover,
\begin{description}
\item[(i)]
For any $t> 0$, the form $\eta_t$ is smooth,
and as $t\arrow 0$, $\eta_t$ converges to $\eta$;
uniformly if $\eta$ is continuous.  

\item[(ii)] The solution $\eta_t$ can be obtained by an integration
with a {\bf heat kernel}, $K_t$, which is an $L^1$-integrable 
$n$-form on $M \times M$,
\[
\eta_t = (\pi_2)_*(\pi_1^* \eta\wedge K_t),
\]
where $\pi_1, \pi_2:\; M \times M \arrow M$ are natural
projections, $\pi_1^*$ the pullback map, and $(\pi_2)_*$
the fiberwise integration. 

\item[(iii)] The form $K_t$ is smooth, for all $t>0$.

\item[(iv)] If $\eta$ is smooth on a neighbourhood of a compact
set $K\subset M$, then all derivatives of $\eta_t$ uniformly
converge to the respective derivatives of $\eta$ on $K$.

\end{description}

\hfill

{\bf Proof:} See e.g. \cite{_Berline_Getzler_Verne_}. \endproof

\hfill

On a K\"ahler manifold, the heat kernel can be locally
expressed as follows.

\hfill

\theorem\label{_heat_ker_Kahl_Theorem_}
Let $M$ be a compact K\"ahler manifold, 
of complex dimension $n$, and $K_t\in \Lambda^{2n}(M)$
the heat kernel. Then
\begin{description}
\item[(i)] $K_t$ has Hodge type $(n,n)$
\item[(ii)] $K_t$ can be decomposed as $K_t = K_t^\circ
  +K_t^{\circ\circ}$, with $K_t^{\circ\circ}$
a smooth form which converges to 0 uniformly
as $t\arrow 0$, and $K^\circ_t$ positive and equal
0 outside of a neighbourhood $U$ of $\diag\subset M\times M$.
Moreover, $U$ can be chosen arbitrarily small.

\item[(iii)] Denote by $H_t$ the operator which maps
$\eta$ to the solution $\eta_t$ of the heat equation.
Then $H_t$ commutes with the operators $d, d^c, \6,
\bar\6$, and with
the multiplication by the K\"ahler form. 

\end{description}

{\bf Proof:} See \cite{_Berline_Getzler_Verne_},
\cite{_Hain:heat_kernel_}. \endproof

\hfill

We are going to approximate 
a $\omega^q$-plurisubharmonic function $\phi$
by $\phi_t:= H_t(\phi)$. To prove \ref{_Heat_Riemannian_Theorem_}, 
we need to show
that $\phi_t$ is also $\omega^q$--plu\-ri\-sub\-har\-monic.
By \ref{_heat_ker_Kahl_Theorem_} (iii), 
it is sufficient to show that the form 
$\eta_t$ is positive for $\eta$ strictly
positive, and $t$ sufficiently small.

\hfill

Using this description of the heat kernel, the following
proposition can be easily established.

\hfill

\proposition \label{_heat_psh_Proposition_}
Let $K\subset M$ be a compact subset of a compact
K\"ahler manifold, and $\phi:\; M \arrow \R$ an
$L^1$-integrable function which is strictly
$\omega^q$-plurisubharmonic on some open 
neighbourhood of $K$. Consider the solutuion
of a heat equation 
\begin{equation}\label{_heat_again_Equation_}
\frac{d\phi_t}{dt}= -\Delta \phi_t, \ \  \phi_0 = \phi,
\end{equation}
Then, for $t$ sufficiently small, $\phi_t$
is $\omega^q$-plurisubharmonic on some open 
neighbourhood of $K$.

\hfill

{\bf Proof:} As we have already mentioned,
to prove \ref{_heat_psh_Proposition_}
it suffices to show that the solution
$\eta_t$  of heat equation is positive, for any
current $\eta\in \Lambda^{q,q}(M, \R)$ which is strictly positive in an open 
neighbourhood $U \supset K$, and $t$ sufficiently small. Replacing
$U$ by a smaller neighbourhood of $K$, we can
always assume that there exists $\epsilon >0$ 
such that $\eta\restrict U > \epsilon \omega^q$.
Writing $K_t = K_t^\circ
  +K_t^{\circ\circ}$ as in \ref{_heat_ker_Kahl_Theorem_},
we can express $\eta_t$ as
\[
\eta_t= \eta_t^\circ + \eta_t^{\circ\circ},
\]
where 
\[
\eta_t^\circ = (\pi_2)_*(\pi_1^* \eta\wedge K_t^\circ),\ \ 
\eta_t^{\circ\circ} = (\pi_2)_*(\pi_1^* \eta\wedge K_t^{\circ\circ}).
\]
Since $K^\circ_t$ is positive, and vanishes outside
of a sufficiently small neighbourhood of $\diag$, the form
$\eta^\circ_t$ is also positive, in some smaller
open neighbourhood $U'\supset K$. Moreover, 
$\eta ^\circ_t > \epsilon \omega^q$ in $U'$
because $\eta > \epsilon \omega^q$, and $H_t(\omega^q)=\omega^q$.
Since $K_t^{\circ\circ}$ is smooth and uniformly
converges to zero,  $\eta_t^{\circ\circ}$ is smooth and 
uniformly converges to zero as $t\arrow 0$. 
Then
\[
\eta_t > \epsilon \omega^q+ \eta_t^{\circ\circ}
\]
is positive for any $t$ for which $|\eta_t^{\circ\circ}|_{L^\infty} < \epsilon$.
For any such $t$, $\phi_t$ is also
$\omega^q$-plurisubharmonic in a neighbourhood of $K$.
\endproof

\hfill

To finish the proof of the local approximation theorem
(\ref{_local_appro_Theorem_}), we shall need the following
lemma.

\hfill

\lemma\label{_open_iso_embe_Kahler_Lemma_}
Let $M$ be an open ball in $\C^n$, centered in 0 and 
equipped with a K\"ahler form $\omega = dd^c \phi$, 
with $\phi$ bounded on $M$. Fix a smaller ball $B
\subsetneq M$, also centered in 0. Then $M$
admits a holomorphic embedding to a compact
K\"ahler manifold $(\tilde M, \tilde \omega)$,
such that $\tilde \omega\restrict B = \omega\restrict B$.

\hfill

{\bf Proof:} In \cite{_Demailly_1982_}, the notion
of a {\em regularized maximum} of two functions was
defined. Choose $\epsilon >0$, and let 
$\max_\epsilon:\; \R^2\arrow \R$ be a smooth, convex
function which satisfies $\max_\epsilon(x, y) = \max(x,y)$
whenever $|x-y|>\epsilon$. Then $\max_\epsilon$ is called
{\bf a regularized maximum}. It is easy to show
(\cite{_Demailly_1982_}) that a regularized maximum
of two strictly plurisubharmonic functions is
again strictly plurisubharmonic.

Replacing the potential $\phi$ by $\phi + \inf_M \phi$,
we can make sure that $0 \leq \phi\leq C$, for some $C\in \R$.
Rescaling, we may assume that the radius of the ball $B$
is 1, and radius of $M$ is $R>1$. Consider the function 
\[  \tilde \phi:=
    \max_\epsilon(\phi, V),
\]
on $M$, where $V(x)= A|x|^2 - A -\epsilon$, where $A$
is a positive number, chosen in such a way
$A R^2 - A - \epsilon > C +\epsilon$. Then
$V(x)< 0$ on $B$, and $V(x)> \phi +\epsilon$ around the
boundary of $M$. The K\"ahler metric 
$\tilde\omega = dd^c \tilde \phi$ is flat 
around the boundary of $M$, and equal $\phi$ on $B$. 
Gluing $(M, \tilde \omega)$ to a round hole in a complex
torus $\tilde M$ with flat K\"ahler metric, we obtain
an isometric embedding 
$(M,\tilde \omega) \hookrightarrow (\tilde M, \tilde \omega)$,
with $\tilde \omega\restrict B = \omega\restrict B$.
This proves \ref{_open_iso_embe_Kahler_Lemma_}. \endproof

\hfill

Now we can finish the proof of the local
approximation theorem (\ref{_local_appro_Theorem_}).
First, the function $\phi$ can be approximated by a
sequence $\phi + \epsilon |z|^2$ of strictly
$\omega^q$-plurisubharmonic functions. Therefore,
it suffices to prove \ref{_local_appro_Theorem_}
assuming that $\phi$ is strictly
$\omega^q$-plurisubharmonic.
Using \ref{_open_iso_embe_Kahler_Lemma_},
we can embed $M$ into a compact torus $\tilde M$,
in such a way that the K\"ahler form $\tilde \omega$
of $\tilde M$ is equal to $\omega$ on
some neighbourhood of $K$. We extend $\phi$
to a continuous function $\tilde \phi$ 
on $\tilde M$, and use the heat equation
to approximate $\tilde \phi$ by
the solutions $\tilde \phi_t$ of the
heat equation. By \ref{_heat_psh_Proposition_}, for
$t$ sufficiently small,
$\tilde \phi_t$ is also $\omega^q$-plurisubharmonic
on some neghbourhood of $K$. This proves
the Local Approximation Theorem
(\ref{_local_appro_Theorem_}). \endproof

\subsection{Functions which are subharmonic on
  $q$-dimensional complex subvarieties}

Let $M$ be an $n$-dimensional Riemannian manifold. A smooth subharmonic
function $\phi:\; M \arrow \R$ is a function which
satisfies $\Delta\phi \geq 0$. This is equivalent to 
$\int_M \Delta\phi \alpha\geq 0$, for any non-negative
positive volume form $\alpha$ with compact support.
However, by Stokes' theorem
\begin{align*}
\int_M \Delta\phi \cdot \alpha=&\int_M *d*d\phi \cdot
\alpha=\int_M d*d\phi \cdot *\alpha=\\ =& -(-1)^{n-1}\int_M
*d\phi \wedge d*\alpha = -\int_M d\phi \wedge *d*\alpha=\\ =& 
\int_M \phi \cdot d*d*\alpha = \int_M \phi \cdot\Delta\alpha.
\end{align*}
Repeating the construction of Subsection 
\ref{_conti_omega^q_psh_Subsection_}, we generalize the
notion of subharmonicity. 

\hfill

\definition\label{_subharmo_Definition_}
A locally $L^1$-integrable
function $\phi:\; M \arrow \R$ is called {\bf subharmonic}
if and only if $\int \phi \Delta\alpha\geq 0$, for
any smooth volume form $\alpha$ with compact support.

This definition is well known, see e.g. 
\cite{_Green_Wu:approx_}. 

\hfill

\theorem \label{_subharmo_and_omega^q_psh_Theorem_}
Let $(M, \omega)$ be a K\"ahler manifold, and 
$\phi:\; M \arrow \R$ a continuous function. 
Then $\phi$ is $\omega^q$-plurisubharmonic if 
and only if $\phi$ is subharmonic on all 
complex submanifolds of dimension 
$q$.\footnote{Here we consider 
complex subvarieties not necessarily closed in $M$.}

\hfill

{\bf Proof:} This result is well known for
the usual plurisubharmonic functions (which are
$\omega^1$-plurisubharmonic in our terminology),
and its proof for general $q$ is no different.

From \cite{_Harvey_Lawson:Psh_}, Theorem 1.4,
it is known that a smooth function which is plurisubharmonic
with respect to a calibration is subharmonic
on all calibrated submanifolds.
For a calibration $\frac 1 {q!}\omega^q$,
calibrated subvarieties are precisely
complex subvarieties of dimension $q$,
hence a smooth $\omega^q$--plurisubharmonic
function is subharmonic on such subvarieties. 

To prove this result for continuous 
plurisubarmonic functions, we use 
the local approximation theorem (\ref{_local_appro_Theorem_}).

To prove the converse statement, assume that
$\phi$ is a continuous function on $M$
which is subharmonic on all complex 
submanifolds of dimension $q$. 
The statement of \ref{_subharmo_and_omega^q_psh_Theorem_}
is local. Replacing $M$ by a smaller open set
if necessary, we may assume that $M$, as a complex
manifold, is isomorphic to a product,
$M= F \times B$, with $\dim F=q$. Denote by
$\pi:\; M \arrow B$ the natural projection.
Let $\pi^*\Vol_B$ be a volume form on $B$ lifted to $M$,
and $\alpha$ a smooth function on $M$ with compact support. 
The cone of strongly positive $(n-q, n-q)$-forms
is locally generated by the forms of type
$\alpha\pi^*\Vol_B$, taken for different $\alpha$ and
$\pi$. Therefore, to prove positivity of the current
$dd^c \phi\wedge \omega^{q-1}$, it suffices
to show that
\[ \int \phi \wedge \omega^{q-1} \wedge
   \pi^*\Vol B \wedge dd^c\alpha \geq 0.
\]
However, this integral can be expressed using the
Fubini formula, as
\begin{equation}\label{_Fubini_subha_Equation_}
\int \phi \wedge \omega^{q-1} \wedge
   \pi^*\Vol B \wedge dd^c\alpha = \int_B \Vol_B \int_F
   \phi \omega^{q-1} \wedge dd^c \alpha.
\end{equation}
Since $\omega^{q-1} \wedge dd^c \alpha= \Delta\alpha
\Vol_F$, and $\phi$ is subharmonic on the fibers of $\pi$,
the second integral is always positive. Therefore,
\eqref{_Fubini_subha_Equation_} is also positive.
We proved \ref{_subharmo_and_omega^q_psh_Theorem_}.
\endproof

\hfill

\remark \label{_subharmo_bigger_than_harmo_Remark_}
There is another definition of subharmonicity: a function
$\phi$ on $M$ is called subharmonic, if for any harmonic
function $\phi_1$, and any compact convex set 
$B\subset M$, with boundary $\6 B$, 
$\phi\restrict{\6 B} \leq \phi_1\restrict{\6 B}$
implies $\phi\restrict{B} \leq \phi_1\restrict{B}$.
This definition is equivalent to
\ref{_subharmo_Definition_},
as shown in \cite{_Green_Wu:approx_}, Lemma 3.1.

\hfill

\corollary \label{_maximum_psh_Corollary_}
A maximum of two continuous $\omega^q$-plurisubharmonic functions
is again $\omega^q$-plurisubharmonic.

\hfill

{\bf Proof:} For subharmonic functions this is well known
(and trivial if we use the definition mentioned in 
\ref{_subharmo_bigger_than_harmo_Remark_}). Now
\ref{_maximum_psh_Corollary_} is implied immediately 
by \ref{_subharmo_and_omega^q_psh_Theorem_}. \endproof

\hfill

In the proof of \ref{_open_iso_embe_Kahler_Lemma_},
a notion of a regularized maximum was defined, following
Demailly (\cite{_Demailly_1982_}). It is well known that a regularized maximum
of two subharmonic functions is again subharmonic.
Using \ref{_subharmo_and_omega^q_psh_Theorem_}, the
following Corollary is immediately obtained.

\hfill

\corollary \label{_reg_maximum_psh_Corollary_}
A regularized maximum of two continuous $\omega^q$--plu\-ri\-sub\-har\-monic functions
is again $\omega^q$-plurisubharmonic.

\endproof

\subsection{Global approximation for
  $\omega^q$-plurisubharmonic functions}

In \cite{_Green_Wu:approx_}, Green and Wu have proven
several powerful theorems which they used to obtain 
the smooth approximation for a number of important
classes of functions on Riemannian manifolds, including
convex, subharmonic, and Lipschitz functions. Their approach
is axiomatic, based on the formal properties that
are satisfied by these classes of functions. 

It turns out that these formal properties are satisfied
also for the sheaf of strictly $\omega^q$-plurisubharmonic
functions as well.

Let $K \subset M$ be a compact subset of a Riemannian
manifold, and $f:\; M \arrow \R$ a continuous function, which is
smooth in a neighbourhood of $K$. Define
\[ 
   {\goth D}(K, i, f):= \sup_K | \nabla^i f|,
\]
where $\nabla:\; (\Lambda^1M)^{\otimes k} \arrow (\Lambda^1M)^{\otimes k+1}$
is the Levi-Civita connection. Consider the pseudonorm
\[
\| f \|_K:= \sup_K |f| + \sum_{i=1}^\infty \frac 1
   {2^i} \max(1, {\goth D}(K, i, f)),
\]
and the corresponding pseudometric $d_k(f, g):=\|f-g\|_K$.

\hfill

\definition
Let ${\goth C}$ be a subsheaf of a sheaf of continuous
functions. We say that ${\goth C}$ {\bf has
  $C^\infty$-stability property} if for any 
$h \in \Gamma({\goth  C}, U)$, and any compact set $K \subset U$,
there exists $\epsilon >0$, such that for any function $f$
on $U$, smooth in a neighbourhood of $K$, and satisfying
$\| f\|_K <\epsilon$, the sum $h+f$, restricted to
some neighbourhood of $K$, lies in ${\goth C}$.

\hfill

\definition
Let ${\goth C}$ be a subsheaf of a sheaf of continuous
functions. We say that ${\goth C}$ {\bf satisfies the
semilocal approximation property}, if for any
\begin{itemize}
\item $\epsilon >0$
\item open subset $U\subset M$
\item compact subset $K \subset U$
\item compact subset $K_1\subset K$ (possibly, empty)
\item function $f\in \Gamma({\goth C}, U)$, which is
  smooth in a neighbourhood of $K_1$
\end{itemize}
there exists a smooth function 
$F\in \Gamma({\goth C}, U)$, such that
\begin{itemize}
\item $\sup_K |f-F| <\epsilon$
\item $d_{K_1}(F, f) <\epsilon$.
\end{itemize}

\hfill

\definition
Let ${\goth C}$ be a subsheaf of a sheaf of continuous
functions. We say that ${\goth C}$ {\bf has local
  approximation property} if for each point $p\in M$, 
there is an open neighbourhood $U\ni p$ such that the restriction
${\goth C}\restrict U$ has the semilocal approximation property.

\hfill

\definition
Let ${\goth C}$ be a subsheaf of a sheaf of continuous
functions. We say that ${\goth C}$ {\bf has the maximum
  closure property} if for any sections 
$f, g\in \Gamma({\goth C}, U)$, the maximum $\max(f,g)$ also belongs to
${\goth C}$. 

\hfill

The following result was proven in
\cite{_Green_Wu:approx_}

\hfill

\theorem\label{_GW_global_app_Theorem_}
(Global Approximation Theorem;
see Theorem 1.1 and Theorem 4.1 of
\cite{_Green_Wu:approx_}).
Let ${\goth C}$ be a subsheaf of a sheaf of continuous
functions satisfying the $C^\infty$-stability, 
maximum closure, and local approximation property.
Then, for any $f\in \Gamma({\goth C}, M)$, and
any positive function $g$ on $M$, there exists
a smooth function $F\in \Gamma({\goth C}, M)$,
such that $|F-f| < g$. 

\hfill

\theorem \label{_Global_appro_omega_psh_Theorem_}
Let $M$ be a K\"ahler manifold, and $f$ a strictly
$\omega^q$--plu\-ri\-sub\-har\-monic continuous function. Then,
for any positive function $g$ on $M$, there exists
a smooth $\omega^q$-plurisubharmonic
function $F\in \Gamma({\goth C}, M)$,
such that $|F-f| < g$.

\hfill

{\bf Proof:} To prove the Global Approximation Theorem
for $\omega^q$--plu\-ri\-sub\-har\-monic functions, we need only to
check that the sheaf of such functions satisfies the 
properties of Green-Wu theorem (\ref{_GW_global_app_Theorem_}). 
The maximum closure property is clear from 
\ref{_maximum_psh_Corollary_}, and $C^\infty$-stability
is clear from the definition of strict
$\omega^q$-plurisubharmonicity.
The local approximation property follows from the heat
equation argument given in Subsection 
\ref{_local_appro_Subsection_}.
Indeed, the solution of heat equation
$\phi_t$ converges to $\phi$ 
uniformly as $t\arrow 0$. Moreover, if $\phi$ is smooth
in a neighbourhood of a compact set $K_1$, then all 
derivatives of $\phi_t$ converge to respective
derivatives of $\phi$.  We proved the global
approximation theorem for continuous strictly
$\omega^q$-plurisubharmonic functions. \endproof


\section{$\omega^q$-plurisubharmonic functions in a
  neighbourhood of a complex variety of dimension $\leq q-1$
}

\subsection{The main statement}

Given a compact complex subvariety $Z \subset M$,
$\dim Z < q$, in some neighbourood of $Z$
there exists a strictly $q$-convex, exhausting function
(\cite{_Barlet_}, \cite{_Demailly:q-convex_}). 
In the present paper, we prove a similar result for 
$\omega^q$-convex functions. Since all $\omega^q$-convex
functions are also $q$-convex, in K\"ahler setting our
theorem also implies this classical result.

\hfill

Recall that a function $\phi:\; M \arrow [-\infty, c[$
is called {\bf exhausting} if for any $c_1< c$
the set $\phi^{-1}([-\infty, c_1])$ is compact.

\hfill

\theorem\label{_omegaq-psh-in-neigh_Theorem_}
Let $M$ be a K\"ahler manifold, and $Z\subset M$ 
a compact complex subvariety, $\dim_\C Z\leq q-1$.
Then there exists an open neighbourhood 
$V\supset Z$ and a smooth strictly 
$\omega^q$-plurisubharmonic function $\phi:\; V \arrow R$
which is exhausting.

\hfill

The proof of \ref{_omegaq-psh-in-neigh_Theorem_}
roughly follows the argument used by J.-P. Demailly in
\cite{_Demailly:q-convex_}. In Subsection \ref{_psh_near_subva_Subsection_},
we construct a strictly $\omega^q$--plu\-ri\-sub\-har\-monic function
in a neighbourhood of $V$, and in Subsection \ref{_almost_psh_Subsection_},
modify it to be exhausting.

\hfill

\remark
\ref{_omegaq-psh-in-neigh_Theorem_} is trivial when $Z$ is
smooth, and $M$ complete as a metric space. Indeed, for a smooth, closed 
submanifold $Z\subset M$, with $TZ$ containing
no $\Phi$-calibrated subspaces, the distance function
$x \arrow \dist(x, Z)$ is $\Phi$-plurisubharmonic
in a neighbourhood of $Z$. This function is 
exhausting in an open neighbourhood
\[\{ x\in M \ \ |\ \ \dist(x, Z) < \epsilon\},\]
for any $\epsilon >0$, because $M$ is complete.

\subsection[Construction of $\omega^q$-plurisubharmonic functions in a
  neighbourhood of a subvariety]{Construction of
  $\omega^q$-plurisubharmonic functions \\ in a
  neighbourhood of a subvariety}
\label{_psh_near_subva_Subsection_}

We use the following lemma

\hfill

\lemma\label{_omega^q_psh_zero_on_bound_Lemma_}
Let $M$ be a compact Riemannian manifold with boundary $\6M$, and
$Z\subset M$ a closed submanifold with a smooth boundary
$\6 Z= Z \cap \6 M$. Assume that 
$M_0:= M\backslash \6M$ is complex and K\"ahler, and
$Z_0:= Z \cap M_0$ is also complex, $\dim_\C Z_0\leq q-1$.
Assume that in a neighbourhood of each point in $\6Z$,
$M$ can be represented as $Z \times B$, where $B$ is 
an open ball,
and this splitting is compatible with the complex structure.
Then there exists a smooth function $\phi:\; U \arrow \R$,
defined on some  neighbourhood $U\supset Z$, which
is strictly $\omega^q$-plurisubharmonic on 
$U_0= U\cap M_0$. Moreover, $\phi$ can be chosen
in such a way that $\phi$ and all its derivatives
vanish on $\6M \cap U$.

\hfill

{\bf Proof:} Notice that only the last statement of
\ref{_omega^q_psh_zero_on_bound_Lemma_} is non-trivial;
if we did not care about the boundary, we could
just take $\phi(x)= \dist(x, Z)$. Passing to a 
covering of $Z$ by sufficiently small open sets,
we find that \ref{_omega^q_psh_zero_on_bound_Lemma_}
is reduced to the following statement, which is
a special case of \ref{_omega^q_psh_zero_on_bound_Lemma_}.

\hfill

\sublemma\label{_omega^q_psh_local_Sublemma_}
Let $Z$ be a smooth Riemannian manifold with a boundary, 
$B$ an open ball in $\C^n$, centered in 0, and 
$M = Z \times B$ their product, equipped with
some Riemannian structure. Assume that 
$M_0 = (Z\backslash \6Z)\times B\subset M$
is K\"ahler, and $\dim_\C Z\leq q-1$. Consider a smooth function
$\theta_Z:\; Z \arrow \R$, positive everywhere
outside of a boundary, and satisfying
\[ 
  |d\theta_Z| < C_1 \theta_Z, \ \
   |\6\bar\6\theta_Z|<C_2\theta_Z,
\]
for some constants $C_1$, $C_2$. Let 
$\phi:\; M \arrow \R$ map $(z, b)$ to $\theta_Z(z)|b|^2$.
Then $\phi\restrict{(Z\backslash \6Z)\times B_1}$ is strictly
$\omega^q$-plurisubharmonic, for some open
ball $B_1\subset B$, also centered in $0$,
with radius depending on the constants
$C_1, C_2$. 

\hfill

{\bf Proof:} Let $\omega_B$ denote the flat
K\"ahler form $\omega_B:= dd^c |b|^2$ on $B$ lifted to $M$,
and $\theta:\; M \arrow \R$ the function $\theta_Z$ lifted to $M$.
Then
\begin{equation}\label{_local_phi_expre_Equation_}
dd^c \phi(z, b) = \theta\cdot \omega_B + dd^c\theta\cdot |b|^2
-\1 \6\theta \wedge \sum_i x_i d \bar x_i +
\1 \bar\6\theta \wedge \sum_i \bar x_i d x_i,
\end{equation}
where $x_i$ are standard coordinates on $B$.
For $|x_i|$ sufficiently small comparing with
$C_1$, the last two terms are small comparing 
with $\theta \omega_B$. Similarly, when $|b|^2$
is small comparing with $C_2$, the sum
$\theta \omega_B + dd^c\theta |b|^2$
obviously satisfies the inequalities 
\eqref{_omega^q-psh_eigen_Equation_},
hence the form
\[
(\theta \omega_B + dd^c\theta |b|^2)\wedge \omega^{q-1} 
\]
is positive. For such $b$, the function $\phi$
is $\omega^q$-plurisubharmonic. We proved
\ref{_omega^q_psh_local_Sublemma_}. \endproof

\hfill

Now \ref{_omega^q_psh_zero_on_bound_Lemma_}
is finished trivially. Consider a locally
finite covering $W_i$ of $Z$, such that some neighbourhood
of $Z$ is covered by  a union of open sets 
$M_i\cong W_i \times B_i$, 
where $B_i$ are isomorphic to open balls in $\C^n$. Taking a function
$\theta_i:\; W_i \arrow \R$ vanishing on a boundary of $W_i$ 
with all its derivatives and satisfying the conditions
of \ref{_omega^q_psh_local_Sublemma_}, we obtain a
strictly $\omega^q$-plurisubharmonic function
on $M_i$ which vanishes on $\6W_i \times B_i$,
with all its derivatives. Summing up these functions,
we obtain a function $\phi$ which is strictly $\omega^q$-plurisubharmonic 
on some neighbourhood of $Z_0$ and vanishes on the boundary 
of $M$, with all its derivatives. \ref{_omega^q_psh_zero_on_bound_Lemma_}
is proven. \endproof

\hfill

\proposition\label{_omega^q_psh_in_neigh_Proposition_}
Let $Z\subset M$ be a closed complex subvariety of a 
K\"ahler manifold. Then there exists a strictly $\omega^q$-plurisubharmonic 
function on an open neighbourhood $U$ of $Z$ in $M$.

\hfill

{\bf Proof:} For $\dim Z=0$, the distance function
$x\arrow \dist(x, Z)$ is strictly plurisubharmonic,
and consequently $\omega^q$-plurisubharmonic in some
neighbourhood of $Z$. Using induction by $\dim Z$, we may
assume that \ref{_omega^q_psh_in_neigh_Proposition_}
is proven for $Z_1 = \Sing(Z)$ (the singular set 
of $Z$). Let $U_1 \supset Z_1$ be an open neighbourhood of 
$Z_1$, and $\phi_1:\; U_1 \arrow \R$ a strictly 
$\omega^q$-plurisubharmonic function. Consider
a smaller neighbourhood $V_1 \subset U_1$, such
that its closure $\bar V_1$ is compact
and contained in $U_1$.\footnote{This property
is usually denoted as $V_1\Subset U_1$.} Let
$\xi:\; M \mapsto [0,1]$ be a smooth
function which is equal 1 on $V_1$
and 0 on $M\backslash U_1$. Consider a 
neighbourhood $W_1 \Subset V_1$ of $Z_1$
with smooth boundary, and let $\phi_0$ be
a strictly $\omega^q$-plurisubharmonic function
defined in a neighbourhood of $Z \backslash W_1$,
and vanishing on $\6W_1$ with all derivatives
(\ref{_omega^q_psh_zero_on_bound_Lemma_}).
We extend $\phi_0$ over $W_1$ by setting 
$\phi_0\restrict {W_1}=0$. Then, $\phi_0$
and $\xi\phi_1$ are defined in an
open neighbourhood $U$ of $Z$.

\begin{figure}[ht]
\begin{center}
\epsfig{file=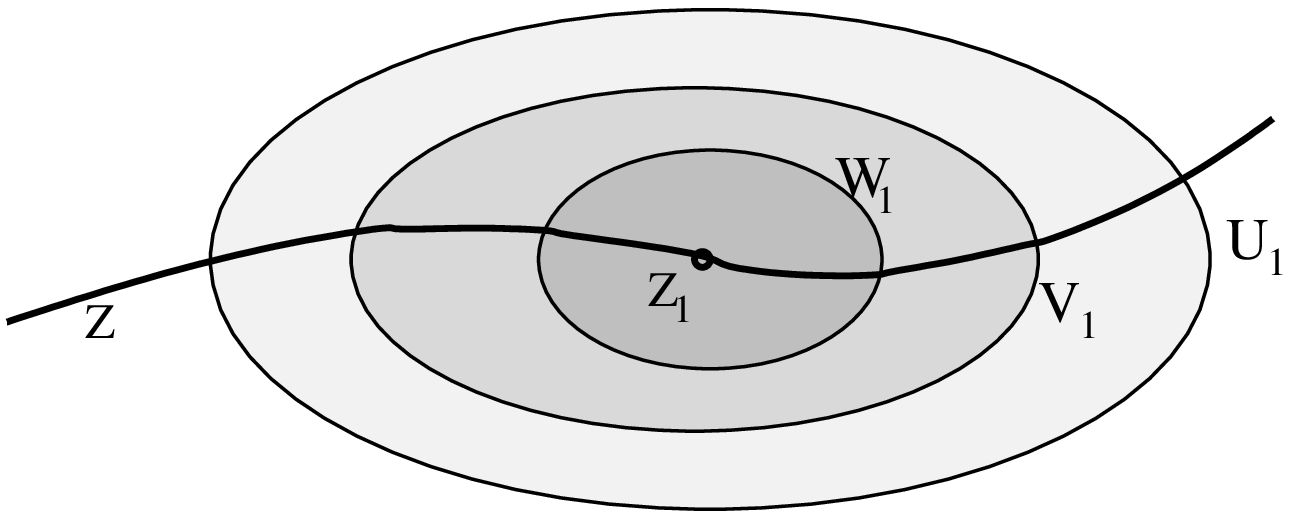,width=0.7\linewidth}
{\parbox[t]{0.7\linewidth}{
\begin{itemize}
\item $V_1$: here $\xi=1$, $\phi_0$ is
  $\omega^q$-plurisubharmonic, and $\phi_1$  strictly 
  $\omega^q$-plurisubharmonic.
\item $U\backslash U_1$: here $\xi=0$, and $\phi_0$ is
strictly  $\omega^q$-plurisubharmonic
\item $U_1\backslash V_1$: here $\phi_0$ is strictly 
  $\omega^q$-plurisubharmonic, hence $C\phi_0 + \xi\phi_1$
is stricty $\omega^q$-plu\-ri\-sub\-har\-monic, for $C\gg 0$.
\end{itemize}}}
\parbox[t]{0.65\linewidth}{{\bf \scriptsize  Figure 1: Gluing
  $\omega^q$-plurisubharmonic functions\\
around the singular set of $Z$}}
\end{center}
\end{figure}

For a sufficiently big constant $C>0$, the sum 
$\phi = C\phi_0 + \xi \phi_1$ is strictly
$\omega^q$-plurisubharmonic. Indeed, on $V_1$,
$\phi= C\phi_0 + \xi \phi_1= C\phi_0 + \phi_1$ 
is  $\omega^q$-plurisubharmonic,
because both of its summands are. Outside of
$U_1$, $\phi= \phi_o$, and this is also 
$\omega^q$-plurisubharmonic. It remains
to show that $\phi$ is $\omega^q$-plurisubharmonic
on the compact set $X:= \overline{U_1\backslash V_1}$,
for $C$ sufficiently big. However, $\phi_0$
is strictly $\omega^q$-plurisubharmonic on this set,
hence satisfies $dd^c \phi_0 \wedge \omega^{q-1}> \epsilon \omega^q$.
Therefore, the sum $C\phi_0 + \xi \phi_1$ 
is  $\omega^q$-plurisubharmonic for 
\[ C> \epsilon^{-1} \|\6\bar\6(\xi \phi_1)\|.\]
Here $\|\6\bar\6(\xi \phi_1)\|(x)=- \min \{\alpha_i\}$,
where $\{\alpha_i\}$ are eigenvalues of the pseudo-Hermitian
form $\6\bar\6\xi \phi_1$ at $x\in X$. We proved
\ref{_omega^q_psh_in_neigh_Proposition_}. \endproof

\subsection{Almost plurisubharmonic functions}
\label{_almost_psh_Subsection_}

Recall that an $L^1$-integrable function 
$\chi:\; M \arrow \R$ on a complex manifold is called
{\bf almost plurisubharmonic} if 
the current $dd^c \chi+ C\omega$ is positive,
for some $C>0$. In \cite{_Demailly_1982_}, Demailly proved the following 
important theorem.

\hfill

\theorem\label{_almost_psh_w_poles_Theorem_}
Let $Z\subset M$ be a compact complex subvariety
of a complex manifold $M$. Then there exists an almost
plurisubharmonic function \[ \chi:\; M \arrow [-\infty, C[\] which
is smooth outside of $Z$ and satisfies 
$\lim\limits_{x\arrow Z}\chi(z) = -\infty$
\endproof

\hfill

Using \ref{_almost_psh_w_poles_Theorem_}, we can finish the proof
of \ref{_omegaq-psh-in-neigh_Theorem_}. This argument is
lifted from \cite{_Demailly:q-convex_}, almost without
changes. Let $\phi:\; V \arrow \R$ be a strictly 
$\omega^q$-plurisubharmonic function in a neighbourhood
of $Z\subset M$. Such a function exists by 
\ref{_omega^q_psh_in_neigh_Proposition_}. 
Rescaling and replacing $V$ with a smaller neighbourhood
if necessary, we may assume that $dd^c\phi \wedge
\omega^{q-1} > \omega^q$. 
Using \ref{_almost_psh_w_poles_Theorem_},
we construct an almost plurisubharmonic
function $\chi$ on $M$, with a pole in $Z$. 
Since $\chi$ is almost plurisubharmonic, 
$dd^c \chi+ C\omega$ is positive. 
Then 
\[ (Cdd^c \phi +dd^c \chi) \wedge \omega^{q-1} >  
C\omega^q +dd^c \chi \wedge \omega^{q-1} = 
(C\omega +dd^c \chi) \wedge \omega^{q-1}.
\]
The last form is positive, because 
$dd^c \chi+ C\omega$ is positive. Therefore,
$\chi_1:=C\phi + \chi$ is $\omega^q$-plurisubharmonic.
Since 
\[ \lim\limits_{x\arrow Z} \chi(x) = -\infty,\]
and $\chi$ is continuous outside of $Z$, for sufficiently
small set $A$, the preimage $W := \chi_1^{-1}([-\infty, A])$ 
contains an open neighbourhood of $Z$.
Assume that $\phi < B-\epsilon$ on $W$. Then 
\[
\psi:= \max_{\frac \epsilon 3} (\phi - B, \chi_1 - A)
\]
is equal to $\chi_1 + A$ in a neighbourhood of $\6 W$, and
to $\phi - B$ in a neighbourhood of $Z$.
Here $\max_{\frac \epsilon 3}$ denotes the regularized
maximum (see the proof of \ref{_open_iso_embe_Kahler_Lemma_}).
By \ref{_reg_maximum_psh_Corollary_}, $\psi$ is strictly
$\omega^q$-plurisubharmonic. The function $\psi$ 
is smooth on outside of $Z$, being a regularized 
maximum of two smooth functions, and it is smooth
in a neighbourhood of $Z$, because 
$\psi = \phi - B$ in a small neighbourhood of $Z$.
Also, $\psi$ is non-positive on $W$, and equal 0 
on $\6W$, hence it is proper and exhausting on $W$. We proved 
\ref{_almost_psh_w_poles_Theorem_}. \endproof

\hfill

We also proved the following result.

\hfill

\corollary\label{_omega_psh_w_pole_Corollary_}
Let $M$ be a K\"ahler manifold, and $Z\subset M$ 
a complex subvariety, $\dim_\C Z\leq q$.
Then there exists an open neighbourhood 
$V$ of $Z$ and a strictly 
$\omega^q$-plurisubharmonic function 
$\chi_1:\; V \arrow [-\infty, 0]$ which is 
smooth outside of $Z$, and satisfies
\[
\lim_{x\arrow Z} \phi(z)=-\infty
\]

\hfill

{\bf Proof:} The function $\chi_1:=C\phi + \chi$ 
defined above obviously satisfies these properties.
\endproof

\hfill


\section{A proof of Sibony's lemma on $L^1$-integrability of 
positive currents}


As an application of \ref{_almost_psh_w_poles_Theorem_},
we give a proof of the following classical result, 
a more general version of which is due to 
Sibony (\cite{_Sibony_}).

\hfill

\theorem\label{_Sibony_Lemma_Theorem_}
(Sibony's lemma). Let $M$ be a K\"ahler manifold,
$\dim_\C M=n$, and $Z\subset M$ a compact
complex subvariety of dimension at most $p-1$.
Consider a weakly  positive, closed, locally $L^1$-integrable
$(n-p, n-p)$-form $\theta$ on $M\backslash Z$. Then $\theta$ 
is locally integrable on $M$.

\hfill

\remark
The usual proof of this lemma does not require $M$
to be K\"ahler, and $Z$ to be compact. 
Instead, one introduces coordinates
and uses the slicing method. In other geometric situations
(on hyperk\"ahler manifolds, for example) there are no
flat coordinates, and slicing fails. Indeed, a
typical hyperk\"ahler manifold does not have
any non-trivial hyperk\"ahler subvarieties, even locally.
The arguments used below for a coordinate-free 
proof of Sibony's lemma can be employed in hyperk\"ahler
geometry, to obtain results on stability of 
derived direct images of coherent sheaves (see
\cite{_V:reflexive_}, \cite{_Verbitsky:Skoda_}).

\hfill

To prove  \ref{_Sibony_Lemma_Theorem_},
we use the following proposition.

\hfill

\proposition\label{_forms_zero_in_neighb_Proposition_}
Let $M$ be a K\"ahler manifold, and $Z\subset M$ a complex
subvariety, $\dim_\C Z<p$. Then there exists an open neighbourhood
$U$ of $Z$, and a sequence $\{\eta_i\}$ of strongly
positive, exact, smooth $(p,p)$-forms on $U$ satisfying
the following.
\begin{description}
\item[(i)] For any open subset $V\subset U$,
with the closure $\bar V$ compact and not 
intersecting $Z$, the restriction 
$\eta_i\restrict V$ stabilizes
as $i\arrow \infty$. Moreover, $\eta_i\restrict V$
is strictly positive for $i\gg 0$ (that is, lies
in the inner part of the strongly positive cone;
see Subsection \ref{_strict_posi_Subsection_} for details).
\item[(ii)] For all $i$, $\eta_i=0$ in some neighbourhood
  of $Z$.
\item[(iii)] The limit $\eta=\lim\eta_i$ is a strictly
positive current on $U$.
\item[(iv)] The forms $\eta_i$ can be written
as $\eta_i = dd^c \phi_i \wedge \omega^{p-1}$,
where $\phi_i$ are smooth functions on $U$.
On any compact set not intersecting $Z$,
the sequence $\{\phi_i\}$ stabilizes as  $i\arrow \infty$.
\end{description}

{\bf Proof:} By \ref{_omega_psh_w_pole_Corollary_},
there exists a strictly $\omega^q$-plurisubharmonic
function $\phi$ on a neighbourhood $U$ of $Z$,
smooth outside of $Z$, and satisfying 
\[
\lim_{x\arrow Z} \phi(z)=-\infty.
\]
Using the regularized maximum, we obtain a sequence
of smooth $\omega^q$-\-plu\-ri\-sub\-har\-monic
functions $\phi_N:= \max_\epsilon(-N, \phi)$. Let
$\eta_i:= dd^c \phi_i \wedge \omega^{p-1}$.
This form is positive, because $\phi_i$ is 
$\omega^q$-plurisubharmonic, and equal to $\phi$
on the set
\[\{ x\in U \ \ | \ \  \phi(x) > -i + \epsilon\}.
\]
This gives the condition (i) of 
\ref{_forms_zero_in_neighb_Proposition_}.
The condition (ii) is also clear, because 
$dd^c\phi_i =0$ on the set
\[
  \{ x\in U \ \ | \ \  \phi(x) < -i - \epsilon\}.
\]
The conditions (iii) and (iv) are apparent from the construction.
We proved \ref{_forms_zero_in_neighb_Proposition_}.
\endproof

\hfill

To finish the proof of \ref{_Sibony_Lemma_Theorem_},
we choose a neighbourhood $U \supset Z$, with 
compact closure and a smooth boundary
$\6U$, and admitting a sequence $\{\eta_i\}$
of smooth forms satisfying \ref{_forms_zero_in_neighb_Proposition_}.
Denote by $\eta$ the current 
$\lim\limits_{i\arrow \infty} \eta_i$.
Since $\eta$ is strictly positive, and $\theta$ is positive,
to prove that $\theta$ 
is locally integrable, it suffices to establish
the inequality 
\begin{equation}\label{_bound_on_inte_Equation_} 
  \int_K \theta\wedge \eta <C
\end{equation}
for any compact set $K\subset U$ not intersecting $Z$,
and a constant $C$ independent from the choice of $K$.
However,
\[
\int_K \theta\wedge \eta= 
\lim\limits_{i\arrow \infty} \int_K \theta\wedge\eta_i,
\]
hence \eqref{_bound_on_inte_Equation_}  would follow
from a similar universal bound $\int_K \theta\wedge \eta_i <C$.
The integral $\int_U \theta\wedge \eta_i <C$ is well defined,
because $\eta_i$ vanish in a neighbourhood of $Z$, and satisfies
\[
\int_K \theta\wedge\eta_i \leq \int_U \theta\wedge \eta_i
\]
for any $K\subset U$. Therefore, to prove  \eqref{_bound_on_inte_Equation_}
it suffices to show that there exists a universal bound
\begin{equation}\label{_int_on_U_eta_i_bounded_Equation_}
\int_U \theta\wedge \eta_i <C,
\end{equation}
with $C$ independent from the choice of $i$.

However,
\begin{equation}\label{_final_int_estimate_Equation_}
\int_U \theta\wedge \eta_i= 
\int_{\6U} \theta\wedge d^c\phi_i \wedge \omega^{p-1}.
\end{equation}
Since $\{\phi_i\}$ stabilizes in a neighbourhood of $\6U$,
as follows from \ref{_forms_zero_in_neighb_Proposition_}
(iv), the integral \eqref{_final_int_estimate_Equation_}
also stabilizes. Therefore, it is universally bounded. 
This proves \eqref{_int_on_U_eta_i_bounded_Equation_}.
We have shown that $\theta$ is integrable.
\endproof

\hfill 

{\bf Acknowledgements:}
I am grateful for F. Reese Harvey and 
H. Blaine Lawson for an interesting e-mail exchange,
and for sending me the Reese Harvey's notes on 
$\omega^q$-plurisubharmonic functions.
Many thanks to Semyon Alesker for his interest in 
\cite{_V:reflexive_}. My gratitude to Semyon Alesker,
Geo Grantcharov and Liviu Ornea for interesting
discussions, to Alexander Rashkovsii for a reference
to \cite{_Khusanov_1_}-\cite{_Khusanov_3_}, and
to Mihnea Col\c toiu for a reference to 
\cite{_Wu:q_complete_} and \cite{_Napier_Ramachandran_}.

\hfill

\hfill

\small{

\noindent {\sc Misha Verbitsky\\
{\sc  Institute of Theoretical and
Experimental Physics \\
B. Cheremushkinskaya, 25, Moscow, 117259, Russia }\\
\tt verbit@mccme.ru }

}

\end{document}